\documentclass[11pt]{amsart}
\usepackage{mathrsfs}
\usepackage{amsfonts}
\usepackage{amsmath}
\usepackage{amssymb}
\usepackage{amsthm}
\usepackage{enumerate}
\usepackage{color}
\usepackage{geometry}
\usepackage{hyperref}
\usepackage{float}
\usepackage{graphicx}
\usepackage{multirow}
\usepackage[numbers,sort&compress]{natbib}
\usepackage{color}
\usepackage[all]{xy}
\usepackage{cases}

\allowdisplaybreaks
\numberwithin{equation}{section}

\theoremstyle{plain}
\newtheorem{prop}{Proposition}[section]

\newtheorem{lemm}[prop]{Lemma}

\newtheorem{theorem}[prop]{Theorem}
\theoremstyle{definition}
\newtheorem{definition}[prop]{Definition}

\newtheorem{example}[prop]{Example}
\newtheorem{remark}[prop]{Remark}

\newcounter{ITEM}
\newcommand\ITEM[1]{\setcounter{ITEM}{#1}\leavevmode\hbox{\rm(\roman{ITEM})}}

\newcommand\PJ{\Pi_{i\in I}^{\ast}J_i}
\newcommand\HS[1]{\leavevmode\null\hspace{#1mm}}
\newcommand\wdots{, ...\HS{0.2}, }
\title{Free Products of digroups}

\author{Guangliang Zhang}
\address{G.Z., School of Mathematics and Systems Science, Guangdong Polytechnic Normal University, Guangzhou 510631, P. R. China}
\email{zgl541@163.com}

\author{Yuqun Chen$^*$}
\address{Y.C., School of Mathematical Sciences, South China Normal University, Guangzhou 510631, P. R. China}
\email{yqchen@scnu.edu.cn}

\author{Zerui Zhang$^{\sharp}$}
\address{Z.Z., School of Mathematical Sciences, South China Normal University, Guangzhou 510631, P. R. China}
\email{zeruizhang@scnu.edu.cn}

\thanks{${}^*$ Supported by the NNSF of China (11571121, 12071156)}

\thanks{${}^{\sharp}$  Supported by Young Teacher Research and Cultivation  Foundation of South China Normal University 20KJ02}

\thanks{${}^{\sharp}$ Corresponding author}



\begin{document}
\begin{abstract}
 We construct the free products of arbitrary digroups,  and thus we solve an open problem of  Zhuchok.
\end{abstract}
\maketitle
\section{Introduction}
A digroup is a generalization of a group, and it first appeared in Loday's work~\cite{Lo99}.  Digroups play an important role in the theory of Leibniz algebras in the problem on finding an appropriate generalization of Lie's third theorem which associates a (local) Lie group to any Lie algebra.  Kinyon made a good progress in terms of digroups to build a partial solution to the above open problem on Leibniz algebras~\cite[Corollary 5.7]{Kinyon04}.

 The study of algebraic properties of digroups has attracted considerable attention.
A digroup is a set with two binary associative operations satisfying some additional
conditions (Definition~\ref{defi-dig}).   Kinyon~\cite{Kinyon04} showed that every digroup is a product of a group and a trivial digroup~\cite{Kinyon04}.
In \cite{Ph}, Phillips offered an equivalent definition for the variety of digroups.
Salazar-D\'\i az, Vel\'asquez and Wills-Toro \cite{Salazar} studied a further generalization of the digroup structure and showed an analogous to the first isomorphism theorem. Ongay, Vel\'asquez and Wills-Toro \cite{Ongay} discussed the notions of normal subdigroups and quotient digroups,
and then  established the corresponding analogues of the classical Isomorphism Theorems. In~\cite{AYZhuchok17},  A. V.  Zhuchok and Y.~V. Zhuchok constructed some interesting digroups. For further investigation on digroups, we refer to
Zhuchok \cite{Zhu17}, Felipe~\cite{Felipe} and Liu~\cite{Liu}.    One of the open problems raised by Zhuchok in \cite{Zhu17} is how to construct the free product of digroups. And our aim of this article is to construct  the free products of arbitrary digroups.

\section{Free products of arbitrary digroups}

Recall that a \textit{disemigroup} $(D, \vdash, \dashv)$ is a set $D$ equipped with two binary operations $\vdash$ and $\dashv$  such that
 the following conditions are satisfied:
\begin{equation}\label{eq00}
\begin{cases}
 \ (D, \vdash)\  \mbox{and} \  (D, \dashv) \  \mbox{are both semigroups}, \\
 \ a\dashv(b\vdash c)=a\dashv (b\dashv c),\\
\ (a\dashv b)\vdash c=(a\vdash b)\vdash c, \\
\ a\vdash(b\dashv c)=(a\vdash b)\dashv c,
\end{cases}
\end{equation}
for all~$a,b,c\in D$. An element $e$ in $D$ is called a \emph{bar-unit} if $e\vdash a = a\dashv e = a$ for every~$a\in D$. And a disemigroup~$(D, \vdash, \dashv)$ is called a dimonoid if~$(D, \vdash, \dashv)$ has a bar-unit. From now on, we shall simply say that $D$ is a disemigroup if $(D, \vdash, \dashv)$ is.

Now we recall an example of a dimonoid, where each element of the dimonoid can be a bar-unit. Moreover, it is in fact a trivial digroup in the sense of \cite{Kinyon04}.
\begin{example}\cite[Example 1.3.b]{Lo99} \label{exam-d-bar}
Let $D=\{e_i\mid i\in I\}$, where~$I$ is an index set. Define~$e_i\vdash e_j=e_j$ and~$e_j \dashv e_i=e_j$ for all~$i,j\in I$. Then~$(D,\vdash,\dashv)$ is a disemigroup and every~$e_i$ can be a bar-unit of~$D$.
\end{example}

\begin{definition}(\cite{Kinyon04, Felipe, Felipe1,Liu,Ph,AYZhuchok17})\label{defi-dig}
 A dimonoid~$D$ with a fixed bar-unit~$e$ is called a \emph{digroup} with respect to~$e$, denoted by~$(D,e)$, if $D$ satisfies the following condition:
$$
 (\forall \ a\in D, \ \exists \  a^{-1}\in D) \ \ \ a \vdash a^{-1} = e = a^{-1} \dashv a.
$$
In this case, we call~$e$ the \emph{fixed bar-unit} of the digroup~$(D,e)$, and call~$a^{-1}$ the \emph{inverse} of~$a$ (with respect to~$e$).
\end{definition}

A more precise notation for the inverse might be~$a^{-1}_{e}$, but since the fixed bar-unit is always clear whenever we write  a digroup~$(D,e)$,  there is no harm to use the notation~$a^{-1}$ instead of~$a^{-1}_{e}$ for the inverse of~$a$.  We refer to~\cite{Felipe1} for more interesting properties of digroups. For instance, from~\cite{Felipe1} we easily deduce that the inverse of~$a$ is unique. We shortly repeat the proof for the convenience of the readers.  Suppose that~$a$ and $b$ are two inverses of~$c$, then we have
$$a=e\vdash a=(b\dashv c)\vdash a=b\vdash (c\vdash a)=b\vdash e
=b\vdash (c\vdash b)=(b\dashv c)\vdash b=e\vdash b=b.$$

Now we recall a useful result which describes the ``group part" of a digroup.  For the convenience of the readers, we quickly repeat the proof.
\begin{lemm}\emph{(\cite{Kinyon04})}\label{Kinyon}
Let $(D,e)$ be a digroup, and let $J=\{a^{-1}\mid a\in D\}$ be the set of all inverses in $D$. Then we have $J=\{a\vdash e \mid a\in D\}$. Moreover, $J$ is a group with unit~$e$ such that the operations~$\vdash$ and~$\dashv$ coincide in $J$.
\end{lemm}
\begin{proof}
 For every~$a\in D$, we have $a^{-1}\vdash (a\vdash e)=(a^{-1}\dashv a)\vdash e=e\vdash e=e$ and
\begin{align*}
(a\vdash e)\dashv a^{-1}&= (a\vdash (a^{-1}\dashv a))\dashv a^{-1}
=((a\vdash a^{-1})\dashv a)\dashv a^{-1}&\\
&=(e\dashv a)\dashv a^{-1}=e\dashv (a\dashv a^{-1})=e\dashv (a\vdash a^{-1})=e.&
\end{align*}
So we deduce $a\vdash e=(a^{-1})^{-1}\in J$ and $\{a\vdash e \mid a\in D\}\subseteq J$. On the other hand, for every~$a^{-1}\in J$, we have
$$
a^{-1}=e\vdash a^{-1}=(a^{-1}\dashv a)\vdash a^{-1}=a^{-1}\vdash (a\vdash a^{-1})=a^{-1}\vdash e.
$$
Thus we obtain $\{a\vdash e \mid a\in D\}=J$.

Now we show that the operations $\vdash$ and~$\dashv$ coincide in $J$. For all $a\vdash e, b\vdash e\in J$, we have
\begin{align*}
(a\vdash e)\vdash (b\vdash e)&= (a\vdash b)\vdash e=(a\vdash  b)\vdash (b^{-1}\dashv b)  =(a\vdash (b\vdash b^{-1}))\dashv b \\
&=(a\vdash e)\dashv (b\dashv e)=(a\vdash e)\dashv (b\vdash e).
\end{align*}
Clearly, $e=e\vdash e\in J$. By the above reasoning, $J$ is a group with the unit~$e$.
\end{proof}
With the notation of Lemma~\ref{Kinyon}, the set $J$ is called the \textit{group part} of~$(D,e)$ and the set
$$ \{d\in D\mid d\vdash a=a\dashv d=a\ \ \mbox{for every }  a\in D\}$$  is called the   \textit{halo part} of the digroup $(D,e)$.

Note that the fixed bar-unit~$e$ plays an important role when one describes a digroup. For instance, let~$D$ be the disemigroup in Example~\ref{exam-d-bar}. If we let~$e_i$ be the fixed bar-unit,  then~$(D, e_i)$ is a digroup and the inverse of~$e_j$ is~$e_i$ for every~$j\in I$. In particular, the group part of $(D, e_i)$ is $\{e_i\}$.

Now we show that the product of an element  in the group part and an element $b$ not in the group part such that the product ``point to" $b$ does not lie in the group part.
\begin{lemm}\label{dj-closed}
Let $(D,e)$ be a digroup with the group part $J$.
  Then for every $a\in J$ and for every~$b\in D\setminus J$, we have $a\vdash b \in D\setminus J$ and~$b\dashv a\in D\setminus J$.
\end{lemm}
\begin{proof}
   If $a\vdash b\in J$, then we have $a^{-1}\vdash (a\vdash b)\in J$ by Lemma \ref{Kinyon}. So we obtain
$$b=e\vdash b=(a^{-1}\dashv a)\vdash b=a^{-1}\vdash (a\vdash b)\in J,$$ which is a contradiction.
Thus we deduce $a\vdash b\in D\setminus J$. Similarly, we have~$b\dashv a\in D\setminus J$.
 \end{proof}

Our aim in this article is to construct the free products of  arbitrary digroups. So we first recall some necessary notions and recall the definition of the free product of digroups.
A map $\varphi$ from a disemigroup $D_1$ to a disemigroup $D_2$ is called a \textit{disemigroup homomorphism}
if  we have $\varphi(a\vdash b)=\varphi(a)\vdash \varphi(b)$ and $\varphi(a\dashv b)=\varphi(a)\dashv \varphi(b)$ for all~$a,b\in D_1$.
Let $(D_1,e_1)$ and $(D_2,e_2)$ be two digroups. A disemigroup homomorphism $\varphi$ from $D_1$ to $D_2$ is called a \textit{digroup homomorphism} from~$(D_1,e_1)$ to $(D_2,e_2)$ if $\varphi(e_1)=e_2$.
\begin{definition}
Let $\{(D_i,e_i) \mid i\in I\}$ be a family of indexed digroups, and let $(D,e)$ be a digroup such that\\
\ITEM1 there exists an injective digroup homomorphism $\varphi_i: (D_i,e_i)\longrightarrow(D,e)$ for every~$i\in I$;\\
\ITEM2 if $(G,e')$ is a digroup, and if there exists a digroup homomorphism $\psi_i: (D_i,e_i)\rightarrow (G,e')$ for every~$i\in I$, then there
exists a unique digroup homomorphism $\eta: (D,e)\rightarrow (G,e')$ such that the diagram
$$
\xymatrix{
  (D_i,e_i) \ar[rr]^{\varphi_i} \ar[dr]_{\psi_i}
                &  &    (D,e) \ar[dl]^{\eta}    \\
                & (G,e')                 }
$$
commutes for every~$i\in I$.

Then $(D,e)$ is called a \textit{free product} of $\{(D_i,e_i)\mid i\in I\}$.
\end{definition}

Clearly, if a free product of digroups $\{(D_i,e_i)\mid i\in I\}$  exists, then it is unique up to isomorphism.

\begin{remark}\label{rema-dj}
Let $\{(D_i,e_i)\mid i\in I\}$ be a family of digroups, and let~$\{(D_i'',e_i'')\mid i\in I\}$ be a family of digroups such that~$(D_i'',e_i'')$ is isomorphic to~$(D_i,e_i)$ for every $i\in I$. Then clearly the free product of~$\{(D_i,e_i)\mid i\in I\}$  is isomorphic to that of~$\{(D_i'',e_i'')\mid i\in I\}$. So there is no harm to assume that~$\{(D_i,e_i)\mid i\in I\}$ is a family of disjoint digroups. From now on, we assume that $\{(D_i,e_i)\mid i\in I\}$ is a family of disjoint digroups and we  shall identify~$(D_i,e_i)$ with~$((D_i\setminus\{e_i\})\cup \{e\}, e)$ for every~$i\in I$ with the obvious multiplication, where~$e$ is a symbol that is not in~$\cup_{i}D_i$.  Then~$(D_i,e_i)$ and~$((D_i\setminus\{e_i\})\cup \{e\}, e)$ are isomorphic as digroups. We denote by~$(D_i', e)$ the digroup $((D_i\setminus\{e_i\})\cup \{e\}, e)$ and denote by $(J_i,e)$ the group part of $(D_i',e)$. Now we fix the notation~$D$ and~$J$ defined by
$$D=\cup_{i\in I}D_i' \mbox{~and~} J=\cup_{i\in I}J_i.$$ We shall construct the free product of~$\{(D_i', e)\mid i\in I\}$ in Theorem~\ref{nffreeproduct}.
\end{remark}

Let~$J^+$ (resp. $J^*$) be the free semigroup (resp. monoid) generated by~$J$, where~$J$ is defined in Remark~\ref{rema-dj}.  Note that~$J_i$ is a group for every~$i\in I$, and for all distinct~$i,t\in I$, we have~$J_i\cap J_t=\{e\}$. Denote by~$\Pi_{i\in I}^{\ast}J_i$ the free product of the groups~$\{J_i\mid i\in I\}$.  Moreover, by the construction of the free products of groups~\cite[Section 2.9]{Jacobson},  we can consider~$J_i$ as a subgroup of~$\Pi_{i\in I}^{\ast}J_i$ for every~$i\in I$.

Now we introduce two special kinds of words over~$D$, which resembles the elements in the free product of semigroups or groups.
\begin{definition}\label{nota-wd}
For every~$a\in D\setminus\{e\}$, there exists a unique index~$i\in I$ such that~$a\in D_i'$, we refer to~$i$ as the index of~$a$, and denote it by~$\sigma(a)$.
A word~$u$ is called a \emph{good word} over~$D$, if~$u=e$, or if~$u=c_1c_2\cdots c_n$, $n\geq 1$,   each $c_t$ lies in~$D\setminus\{e\}$ and $\sigma(c_p)\neq \sigma(c_{p+1})$ for every integer~$p\leq n-1$.  Denote by~$W(D)$ the set of all good words over~$D$.  A word~$u$ (over~$D$) is called a \emph{reduced word} in~$J^+$ if~$u$ lies in~$W(D)\cap J^+$.
\end{definition}

 We also note that,  the set of all reduced words in~$J^{+}$ forms a set of normal forms of elements in~$\PJ$. And this is why we use the name ``reduced words".
  More precisely, let~$\theta$ be a semigroup homomorphism from~$J^+$ to~$\PJ$ such that~$\theta(a)=a$ for every~$a\in J$. Then clearly~$\theta$ is a surjective homomorphism.  So there exists a congruence~$\rho$ on~$J^+$ such that we have~$J^+/\rho\cong \PJ$. By the construction of the free products of arbitrary groups, for every word~$u\in J^+$, there exists a unique reduced word~$\underline{u}$ in~$J^+$ such that we have~$u\rho=\underline{u}\rho$ in $J^+/\rho$. For instance, let $a\in J_1$, $b\in J_2\setminus\{e\}$ such that~$aa=c\neq e$ in~$J_1$. Then we have~$\underline{aab}=cb$.
In particular, if~$u$ is a reduced word, then we have~$\underline{u}=u$ in~$ J^+$.
Moreover, it follows that~$\underline{\ }$ is a map from~$J^+$ to~$J^+$ such that for all~$u,v\in J^+$, we have
  $$\underline{\underline{u}\ \underline{v}} =\underline{\underline{u}v} =\underline{u \underline{v}} = \underline{uv} \mbox{ and } \underline{eu}=\underline{ue}=\underline{u}. $$

Now we introduce a map~$^{\sharp}$ from~$D^*$ to $D^*$.  For every $a\in D$,   we define
$$
a^{\sharp}=a\vdash e.
$$
Note that, if~$a\in D_i'$, then~$a\vdash e$ means the product in~$D_i'$. Since all the $D_i'\setminus\{e\}$   and~$D_j'\setminus\{e\}$ are disjoint for all~$i\neq j\in I$, the notation $a\vdash e$ makes sense and $a\vdash e$ means an element in~$D_i'$. By Lemma \ref{Kinyon}, we have~$a^{\sharp}=a\vdash e\in J\subseteq J^+$.
Moreover, for every word~$u=a_1a_2...a_m\in D^{+}$, where~$D^+$ is the free semigroup generated by~$D$ and every~$a_t$ lies in~$D$, we define
$$u^{\sharp}=a_1^{\sharp}a_2^{\sharp}...a_m^{\sharp}$$
to be an element in~$J^+\subseteq D^*$. Clearly, $J^+$ is a subsemigroup of~$D^*$, where~$D^*$ is the free monoid generated by~$D$. Then we obtain that
$\underline{u^{\sharp}}=\underline{a_1^{\sharp}a_2^{\sharp}...a_m^{\sharp}}$  is a reduced word  in~$J^+$.  Now we extend~$^{\sharp}$ by defining
$$\varepsilon^{\sharp}=e,$$ where~$\varepsilon$ is the empty word.
Then $^{\sharp}$ becomes a map from~$D^*$ to~$D^*$ such that the image of~$D^*$ is contained in $J^+$.

The following result is an obvious corollary of Lemma~\ref{Kinyon}, and thus the proof is omitted.
\begin{lemm}\label{Jsharp}
The map $^{\sharp}$ is a semigroup homomorphism from $D^{*}$ to $D^{*}$ such that for every~$a\in J$, we have~$a^{\sharp}=a$. Moreover, for every~$u\in D^*$, we have $({u}^{\sharp})^{\sharp}=u^{\sharp}\in J^+$.
\end{lemm}

Recall that for every disemigroup~$M$, for all~$a_1,...,a_t$ in~$M$, every parenthesizing of
$$
a_1\vdash a_2\cdots\vdash a_m \dashv \cdots\dashv a_t
$$
 gives the same element in~$M$~\cite{Lo99}. In light of this, we define a formal expression~$u\dot{a}v$ as follows: for all reduced words~$u,v\in J^+$ and for every~$a\in D$, say~$u=a_1...a_m$ and~$v=b_1...b_n$, where all~$a_1,...,a_m,b_1,...,b_n$ lie in~$J$,  $m\geq 0$ and~$n\geq 0$,  we define
 $$u\dot{a}v=a_1...a_m\dot{a}b_1...b_n,$$
where if~$m$ (resp.~$n$) is 0, that is, when $u$ (resp.~$v$) is the empty word,  it means that $u$ (resp.~$v$) does not appear. In particular, if $u$ and $v$ are both the empty word, then by~$u\dot{a}v$ we mean $\dot{a}$. Therefore, for all words~$u,v\in D^*$ and for every~$a\in D$, the formula
 $\underline{u^{\sharp}}\dot{a}\underline{v^{\sharp}}$ makes sense.

 For the moment, the formula $a_1...a_m\dot{a}b_1...b_n$ is just a formal expression that does not lie in any digroup. But we shall construct a set consisting
of some of such elements and define certain products on the elements of the set to make it a digroup such that
 we have
 $$a_1...a_m\dot{a}b_1...b_n=a_1\vdash \cdots \vdash a_m\vdash a\dashv b_1\dashv \cdots \dashv b_n.$$
Suppose that we already have the free product of~$\{(D_i',e)\mid i\in I\}$. Then for all~$a_i\in D_i'$ and~$a_j\in D_j'$,   in the free product, we have
$$a_i\vdash a_j=(a_i\dashv e)\vdash a_j= (a_i\vdash e)\vdash a_j=a_i^\sharp\vdash a_j$$
and
$$a_i\dashv a_j=a_i\dashv (a_j\dashv e)=a_i\dashv (a_j\vdash e)=a_i\dashv a_j^\sharp.$$
Therefore, for an arbitrary element $a_1\vdash \cdots \vdash a_m\vdash a\dashv b_1\dashv \cdots \dashv b_n$ in the free product of~$\{(D_i',e)\mid i\in I\}$, we may always assume~$a_p\in J$ and~$b_t\in J$ for all~$p\leq m$ and~$t\leq n$. This motivates us to introduce the notion of good center-words as follows. Recall that~$W(D)$ is the set of all good words over~$D$ (Definition~\ref{nota-wd}).
We call $u\dot{a}v$ a \emph{good center-word} over~$D$
  if we have

   \ITEM1 $u,v\in J^*$, $a\in D$, $uav\in W(D)$;

   \ITEM2  if   $a\in J$, then $u$ is the empty word.

Denote by~$CW(D)$ the set of all good center-words over~$D$, that is,
\begin{equation}\label{nota-cwd}
CW(D)=\{u\dot{a}v \mid u,v\in J^*,a\in D, uav\in W(D); \mbox{ if } a\in J, \mbox{ then } u=\varepsilon\}.
\end{equation}
Now we introduce an extended reduction rule on certain ``center-words" as follows.
For all reduced words $u=a_1...a_m, v=b_1...b_n\in J^+$, where all~$a_1,...,a_m,b_1,...,b_n$ lie in~$J$, and for every~$a\in D\setminus J$, we define
$$
[u\dot{a}v]=
\begin{cases}
u\dot{a}v,   &   \mbox{if}~a_m,a~\mbox{do not lie in the same set}~D_i'~\mbox{for any}~i\in I, \mbox{and}~a,b_1\\
            &~\mbox{do not lie in the same set}~D_j'~\mbox{for any}~j\in I;   \\
a_1...a_{m-1}\dot{c}v, &  \mbox{if}~a_m,a\in D_i'~\mbox{for some index}~i\in I,~a_m\vdash a=c\in D_i',~b_1\notin D_i'; \\
u\dot{c}b_2...b_n, &  \mbox{if}~a,b_1\in D_i'~\mbox{for some index}~i\in I,~a\dashv b_1=c\in D_i',~a_m\notin D_i'; \\
 a_1...a_{m-1}\dot{c}b_2...b_n, &  \mbox{if}~a_m,a,b_1\in D_i'~\mbox{for some index}~i\in I, (a_m\vdash a)\dashv b_1=c\in D_i'.
\end{cases}
$$
By Lemma~\ref{dj-closed}, it follows that~$[u\dot{a}v]$ is a good center-word in~$CW(D)$.
Finally, for every reduced word~$u=a_1...a_m \in  J^+$, where~$a_1\wdots a_m$ lie in~$J$,  we define
$$[u]=\dot{a}_1...a_m.$$
In particular, for every word~$u\in D^*$, the formula
 $[\underline{u^{\sharp}}]$ makes sense.

The next lemma shows that the extended reduction rule is compatible with the reduction rule in~$\PJ$ in the following sense.
\begin{lemm}\label{bracket-red}
For all reduced words~$u=a_1...a_m$ and~$v=b_1...b_n\in  J^+$,  where all~$a_i,b_j$ lie in~$J$, $m\geq0$ and~$n\geq 0$, for every~$a\in D\setminus J$, if~$[u\dot{a}v]=c_1...c_p\dot{c}d_1...d_q$, then we have~$
\underline{a_1^{\sharp}...a_m^{\sharp}a^{\sharp}b_1^{\sharp}...b_n^{\sharp}}
=\underline{c_1^{\sharp}...c_p^{\sharp}c^{\sharp}d_1^{\sharp}...d_q^{\sharp}}$.
\end{lemm}
\begin{proof}
For every index~$i\in I$, we have the following three observations:
 \ITEM1 If~$ a_m\vdash a=c$ in~$D_i'$,  then since by Lemma~\ref{Kinyon} the operations~$\vdash$ and~$\dashv$ coincide in~$J_i$, we have
 $$c^{\sharp}=(a_m\vdash a)\vdash e=(a_m\vdash e) \vdash (a\vdash e)=\underline{a_m^{\sharp}a^{\sharp}}. $$
 \ITEM2 If~$ a\dashv b_1=c$ in $D_i'$, then we have
 $$c^{\sharp}=(a\dashv b_1)\vdash e=(a\vdash e) \vdash (b_1\vdash e) =\underline{a^{\sharp}b_1^{\sharp}}. $$
 \ITEM3 If~$(a_m\vdash a)\dashv b_1=c$ in $D_i'$, then we have
 $$c^{\sharp}=((a_m\vdash a)\dashv b_1)\vdash e=(a_m\vdash e)\vdash (a\vdash e)\vdash (b_1\vdash e)=\underline{a_m^{\sharp}a^{\sharp}b_1^{\sharp}}.$$
 The proof of the lemma follows immediately by the definition of~$[u\dot{a}v]$.
\end{proof}

With the above notations, we can define new operations~$\vdash$ and~$\dashv$ on the set of  all good center-words over $D$ in Definition~\ref{defi-fp}. Since the operations extend those of~$D_i'$ for every~$i\in I$, we use the same symbol~$\vdash$ and~$\dashv$ but not~$\vdash',\dashv'$.  We shall prove that the following definition turns out to be a way of  constructing the free products  of arbitrary digroups in Theorem~\ref{nffreeproduct}.
\begin{definition}\label{defi-fp}
With the notion of Remark~\ref{rema-dj}, let~$CW(D)$ be the set of all good center-words over $D$ defined in~\eqref{nota-cwd}. We define two operations~$\vdash$ and~$\dashv$ on~$CW(D)$ as follows. For all~$u_1\dot{a}_1v_1,u_2\dot{a}_2v_2\in CW(D)$, we define
$$u_1\dot{a}_1 v_1\vdash u_2\dot{a}_2 v_2 =
\begin{cases}
[\underline{u_1^{\sharp}a_1^{\sharp}v_1^{\sharp}u_2^{\sharp}}\dot{a}_2v_2^{\sharp}],  & \mbox{if}~a_2\in D\setminus J;\\
&\\
[\underline{u_1^{\sharp}a_1^{\sharp}v_1^{\sharp}u_2^{\sharp}a_2^{\sharp}v_2^{\sharp}}],  & \mbox{if}~a_2\in J,
\end{cases}
$$
and
$$u_1\dot{a}_1 v_1\dashv u_2\dot{a}_2 v_2 =
\begin{cases}
[u_1^{\sharp}\dot{a}_1\underline{v_1^{\sharp}u_2^{\sharp}a_2^{\sharp}v_2^{\sharp}}],  & \mbox{if}~a_1\in D\setminus J;\\
&\\
[\underline{u_1^{\sharp}a_1^{\sharp}v_1^{\sharp}u_2^{\sharp}a_2^{\sharp}v_2^{\sharp}}],  & \mbox{if}~a_1\in J.
\end{cases}
$$
\end{definition}

Before showing that~$(CW(D), \vdash, \dashv)$ is a digroup, we first establish several useful properties of~$\vdash$ and~$\dashv$.  Now we show that,
 when the product ``points to" an element, we may assume that the other element is of certain form.
\begin{lemm}\label{mul-sharp}
With the notation of~Definition~\ref{defi-fp}, for all~$u_1\dot{a}_1v_1,u_2\dot{a}_2v_2\in CW(D)$, we have~$u_1\dot{a}_1v_1\vdash u_2\dot{a}_2v_2=[\underline{u_1^{\sharp}a_1^{\sharp}v_1^{\sharp}}]\vdash u_2\dot{a}_2v_2 $
and $u_1\dot{a}_1v_1\dashv u_2\dot{a}_2v_2=u_1\dot{a}_1v_1\dashv [\underline{u_2^{\sharp}a_2^{\sharp}v_2^{\sharp}}] $.
\end{lemm}
\begin{proof}
We just prove~$u_1\dot{a}_1v_1\vdash u_2\dot{a}_2v_2=[\underline{u_1^{\sharp}a_1^{\sharp}v_1^{\sharp}}]\vdash u_2\dot{a}_2v_2 $ because the other one can be proved similarly. If~$a_2\in J$, then we obtain
  $$u_1\dot{a}_1v_1\vdash u_2\dot{a}_2v_2
  =[\underline{u_1^{\sharp}a_1^{\sharp}v_1^{\sharp} u_2^{\sharp}a_2^{\sharp}v_2^{\sharp}}]
  =[\ \underline{\underline{u_1^{\sharp}a_1^{\sharp}v_1^{\sharp}} u_2^{\sharp}a_2^{\sharp}v_2^{\sharp}}\ ]
  =[\underline{u_1^{\sharp}a_1^{\sharp}v_1^{\sharp}}]\vdash u_2\dot{a}_2v_2.$$
If~$a_2\notin J$, then we have
$$
  u_1\dot{a}_1v_1\vdash u_2\dot{a}_2v_2
   = [\underline{u_1^{\sharp}a_1^{\sharp}v_1^{\sharp} u_2^{\sharp}}\dot{a}_2v_2^{\sharp}]=[\underline{\underline{u_1^{\sharp}a_1^{\sharp}v_1^{\sharp}} u_2^{\sharp}}\dot{a}_2v_2^{\sharp}]
   =[\underline{u_1^{\sharp}a_1^{\sharp}v_1^{\sharp}}]\vdash u_2\dot{a}_2v_2.
$$
The proof is completed.
\end{proof}

The next two lemmas show that the extended reduction rule is compatible with the operations $\vdash$ and $\dashv$.
\begin{lemm}\label{ass1}
For all reduced words~$u,v\in  J^+$, and for every~$a\in D\setminus J$, we have~$[\underline{uv}\dot{a}]
=[u]\vdash[v\dot{a}]$ and
$[\dot{a}\underline{uv}]
=[\dot{a}u]\dashv[v]$. Moreover, for every reduced word~$w\in  J^+$, if~$aw$ is a good word or if~$w=e$, then we have~$[\underline{uv}\dot{a}w]=
[u]\vdash[v\dot{a}w]$; and if~$wa$ is a good word or if~$w=e$, then we have~$[w\dot{a}\underline{uv}]
=[w\dot{a}u]\dashv[v]$.
\end{lemm}
\begin{proof}
We first show~$[\underline{uv}\dot{a}]=[u]\vdash[v\dot{a}]$.
  Assume that~$u=a_1...a_m$, $v=b_1...b_n$, where all~$a_p,b_q$ lie in~$J$ and~$m,n$ are positive integers. There are several cases to consider:

  Case 1: If~$a,b_n\notin D_j'$ for any~$j\in I$, then by Definition~\ref{defi-fp}, we have
  $$[u]\vdash[v\dot{a}]= \dot{a}_1...a_m \vdash v\dot{a}=[\underline{uv}\dot{a}].$$

  Case 2: If~$a,b_n\in D_j' $ for some~$j\in I$ such that~$b_n\vdash a=b$ in~$D_j'$, and if~$\underline{uv}=\underline{ub_1...b_{n-1}}b_n$ as words in the free semigroup~$J^+$,  then by Definition~\ref{defi-fp}, we have
  $$[u]\vdash[v\dot{a}]= \dot{a}_1...a_m \vdash b_1...b_{n-1}\dot{b}=[\underline{ub_1...b_{n-1}}\dot{b}]
  =[\underline{ub_1...b_{n-1}}b_n\dot{a}]=[\underline{uv}\dot{a}].$$

   Case 3: If~$a,b_n\in D_j' $ for some~$j\in I$ such that~$b_n\vdash a=b$ in~$D_j'$, and if~$\underline{uv}\neq \underline{ub_1...b_{n-1}}b_n$ as words in the free semigroup~$J^+$, then we deduce
 $$\underline{a_1...a_{m}b_1...b_{n-1}}=a_1...a_{m-n}a_{m-n+1}$$
 as words in the free semigroup~$J^+$,
 and we have~$a_{m-n+1},b_n,a\in D_j'$.
 Moreover, for every~$0\leq k\leq n-2$, we have~$b_{k+1}=a_{m-k}^{-1}\in J_{i_{k}}$ for some~$i_k\in I$. If~$b_n=a_{m-n+1}^{-1}$, then we have
 $\underline{uv}=a_{1}...a_{m-n}$ and
  $$a_{m-n+1}\vdash b=a_{m-n+1}\vdash (b_n\vdash a)
  =(a_{m-n+1}\vdash b_n)\vdash a=e\vdash a= a.$$ So we deduce
  \begin{align*}
   [u]\vdash[v\dot{a}]=&\dot{a}_1...a_m \vdash b_1...b_{n-1}\dot{b}
 =[\underline{a_1...a_mb_1...b_{n-1}}\dot{b}]=[a_{1}...a_{m-n}a_{m-n+1}\dot{b}]&\\
    =&a_{1}...a_{m-n}\dot{a}=[\underline{uv}\dot{a}].&
 \end{align*}
 If~$b_n\neq a_{m-n+1}^{-1}$, then we have~$a_{m-n+1}\vdash b_n=c\neq e$. Assume~$c\vdash a=d$. Then we have~$\underline{uv}=a_{1}...a_{m-n}c$  in the free semigroup~$J^+$ and
 $$a_{m-n+1}\vdash b=a_{m-n+1}\vdash (b_n\vdash a)=(a_{m-n+1}\vdash b_n)\vdash a
 =c\vdash a=d.$$ So we deduce
 \begin{align*}
   [u]\vdash[v\dot{a}]=&\dot{a}_1...a_m \vdash b_1...b_{n-1}\dot{b}
 =[\underline{a_1...a_mb_1...b_{n-1}}\dot{b}]=[a_{1}...a_{m-n}a_{m-n+1}\dot{b}]&\\
    =&a_{1}...a_{m-n}\dot{d}=[a_{1}...a_{m-n}c\dot{a}]=[\underline{uv}\dot{a}].&
 \end{align*}
 Similarly, we can show~$[\dot{a}\underline{uv}]
=[\dot{a}u]\dashv[v]$.

 As for what remained, if~$w=e$, then we have~$[u\dot{a}w]=[u\dot{a}]$ and~$[w\dot{a}u]=[\dot{a}u]$. So the claim follows by the above proof; if~$w\neq e$, then by Lemma~\ref{dj-closed}, what remained can be proved similarly by noting that~$w$ is never involved in the above reasoning on the extended reduction rule.
\end{proof}

\begin{lemm}\label{ass2}
  For all reduced words~$u,v\in  J^+$, and for every~$a\in D\setminus J$, we have~$[u\dot{a}v]=[u\dot{a}]\dashv [v]=[u]\vdash [\dot{a}v]$.
\end{lemm}
\begin{proof}
By the definition of~$[u\dot{a}v]$ and by Lemma~\ref{dj-closed}, we clearly have~$[u\dot{a}v]=[u\dot{a}]\dashv [v]$.
Similarly, by the fact that~$(b\vdash c)\dashv d=b\vdash (c\dashv d)$ for all~$b,c,d$ in a disemigroup,   we easily deduce~$[u\dot{a}v]=[u]\vdash [\dot{a}v]$.
\end{proof}

Now we show that~$CW(D)$ is a disemigroup with respect to the operations defined in Definition~\ref{defi-fp}.
\begin{lemm}\label{pre-dise}
For all~$u_i\dot{a}_iv_i\in CW(D)$, $1\leq i\leq 3$, we have

\ITEM1~$(u_1\dot{a}_1v_1\vdash u_2\dot{a}_2v_2)\vdash u_3\dot{a}_3v_3= (u_1\dot{a}_1v_1\dashv u_2\dot{a}_2v_2)\vdash u_3\dot{a}_3v_3$;

\ITEM2~$u_1\dot{a}_1v_1\dashv (u_2\dot{a}_2v_2\dashv u_3\dot{a}_3v_3)= u_1\dot{a}_1v_1\dashv (u_2\dot{a}_2v_2\vdash u_3\dot{a}_3v_3)$;

\ITEM3~$(u_1\dot{a}_1v_1\vdash u_2\dot{a}_2v_2)\vdash u_3\dot{a}_3v_3=u_1\dot{a}_1v_1\vdash (u_2\dot{a}_2v_2\vdash u_3\dot{a}_3v_3)$;

\ITEM4~$(u_1\dot{a}_1v_1\dashv u_2\dot{a}_2v_2)\dashv u_3\dot{a}_3v_3= u_1\dot{a}_1v_1\dashv (u_2\dot{a}_2v_2\dashv u_3\dot{a}_3v_3)$;

\ITEM5~$(u_1\dot{a}_1v_1\vdash u_2\dot{a}_2v_2)\dashv u_3\dot{a}_3v_3= u_1\dot{a}_1v_1\vdash (u_2\dot{a}_2v_2\dashv u_3\dot{a}_3v_3)$.\\
In particular, $CW(D)$ is a disemigroup.
\end{lemm}
\begin{proof}
We first prove Point~\ITEM1.
 For every~$\delta\in \{\vdash, \dashv\}$, assume~$(u_1\dot{a}_1v_1) \delta (u_2\dot{a}_2v_2)=u\dot{a}v$.
 Then by Lemma~\ref{bracket-red} and by Definition~\ref{defi-fp}, we have
 $$\underline{u^{\sharp}a^{\sharp}v^{\sharp}}
 =\underline{u_1^{\sharp}a_1^{\sharp}v_1^{\sharp}u_2^{\sharp}a_2^{\sharp}v_2^{\sharp}} \ \mbox{and thus }
 [\underline{u^{\sharp}a^{\sharp}v^{\sharp}}]
 =[\underline{u_1^{\sharp}a_1^{\sharp}v_1^{\sharp}u_2^{\sharp}a_2^{\sharp}v_2^{\sharp}}].$$
So by Lemma~\ref{mul-sharp},
   we deduce
   \begin{align*}
     ((u_1\dot{a}_1v_1) \delta (u_2\dot{a}_2v_2))\vdash u_3\dot{a}_3v_3
   &=u\dot{a}v \vdash u_3\dot{a}_3v_3
   =[\underline{u^{\sharp}a^{\sharp}v^{\sharp}}]\vdash u_3\dot{a}_3v_3&\\
   &=[\underline{u_1^{\sharp}a_1^{\sharp}v_1^{\sharp}u_2^{\sharp}
   a_2^{\sharp}v_2^{\sharp}}]\vdash u_3\dot{a}_3v_3.&
   \end{align*}
   In particular, we have~$(u_1\dot{a}_1v_1\vdash u_2\dot{a}_2v_2)\vdash u_3\dot{a}_3v_3= (u_1\dot{a}_1v_1\dashv u_2\dot{a}_2v_2)\vdash u_3\dot{a}_3v_3$. Point~\ITEM2 can be proved similarly.

   Now we prove Point~\ITEM3.  If~$a_3\notin J$, then by the above reasoning and by Lemmas~\ref{mul-sharp} and~\ref{ass1}, we have
\begin{align*}
  u_1\dot{a}_1v_1\vdash (u_2\dot{a}_2v_2\vdash u_3\dot{a}_3v_3)
=&u_1\dot{a}_1v_1\vdash [\underline{u_2^{\sharp}a_2^{\sharp}v_2^{\sharp}u_3^{\sharp}}\dot{a}_3v_3^{\sharp}]
=[\underline{u_1^{\sharp}a_1^{\sharp}v_1^{\sharp}}]\vdash [\underline{u_2^{\sharp}a_2^{\sharp}v_2^{\sharp}u_3^{\sharp}}\dot{a}_3v_3^{\sharp}]&\\
=&[\underline{u_1^{\sharp}a_1^{\sharp}
v_1^{\sharp}u_2^{\sharp}a_2^{\sharp}v_2^{\sharp}u_3^{\sharp}}\dot{a}_3v_3^{\sharp}] ~~(\mbox{by Lemma~\ref{ass1}})&\\
=&[\underline{u_1^{\sharp}a_1^{\sharp}
v_1^{\sharp}u_2^{\sharp}a_2^{\sharp}v_2^{\sharp}}]\vdash u_3\dot{a}_3v_3
=(u_1\dot{a}_1v_1\vdash u_2\dot{a}_2v_2)\vdash u_3\dot{a}_3v_3.&
\end{align*}
If~$a_3\in J$, then by Lemma~\ref{mul-sharp} and by the above reasoning, we have
\begin{align*}
  u_1\dot{a}_1v_1\vdash (u_2\dot{a}_2v_2\vdash u_3\dot{a}_3v_3)
=&u_1\dot{a}_1v_1\vdash [\underline{u_2^{\sharp}
a_2^{\sharp}v_2^{\sharp}u_3^{\sharp}a_3^{\sharp}v_3^{\sharp}}]&\\
=&[\underline{u_1^{\sharp}a_1^{\sharp}v_1^{\sharp}u_2^{\sharp}
a_2^{\sharp}v_2^{\sharp}u_3^{\sharp}a_3^{\sharp}v_3^{\sharp}}]
=[\underline{u_1^{\sharp}a_1^{\sharp}
v_1^{\sharp}u_2^{\sharp}a_2^{\sharp}v_2^{\sharp}}]\vdash u_3\dot{a}_3v_3&\\
=&(u_1\dot{a}_1v_1\vdash u_2\dot{a}_2v_2)\vdash u_3\dot{a}_3v_3.&
\end{align*}
Point~\ITEM4 can be proved similarly.

Finally we prove Point~\ITEM5. If~$a_2\in J$, then clearly we have
  $$(u_1\dot{a}_1v_1\vdash u_2\dot{a}_2v_2)\dashv u_3\dot{a}_3v_3=
   [\underline{u_1^{\sharp}a_1^{\sharp}v_1^{\sharp}u_2^{\sharp}a_2^{\sharp}v_2^{\sharp}
   u_3^{\sharp}a_3^{\sharp}v_3^{\sharp}}]=u_1\dot{a}_1v_1\vdash (u_2\dot{a}_2v_2\dashv u_3\dot{a}_3v_3).$$
   If~$a_2\notin J$, then by Lemma~\ref{ass2}, we have
   $$u_1\dot{a}_1v_1\vdash u_2\dot{a}_2v_2
   =[\underline{u_1^{\sharp}a_1^{\sharp}v_1^{\sharp}u_2^{\sharp}}\dot{a}_2v_2^{\sharp}]
   =[\underline{u_1^{\sharp}a_1^{\sharp}v_1^{\sharp}u_2^{\sharp}}\dot{a}_2]\dashv [v_2^{\sharp}]$$
   and
   $$u_2\dot{a}_2v_2\dashv u_3\dot{a}_3v_3
   =[u_2^{\sharp}\dot{a}_2\underline{v_2^{\sharp}u_3^{\sharp}a_3^{\sharp}v_3^{\sharp}}]
   =[u_2^{\sharp}]\vdash [\dot{a}_2\underline{v_2^{\sharp}u_3^{\sharp}a_3^{\sharp}v_3^{\sharp}}]. $$
 By Points~\ITEM3 and~\ITEM4, and by Lemma~\ref{ass2} again, we deduce
\begin{align*}
 &(u_1\dot{a}_1v_1\vdash u_2\dot{a}_2v_2)\dashv u_3\dot{a}_3v_3 =([\underline{u_1^{\sharp}a_1^{\sharp}v_1^{\sharp}u_2^{\sharp}}\dot{a}_2]
   \dashv [v_2^{\sharp}])\dashv u_3\dot{a}_3v_3&\\
    =&[\underline{u_1^{\sharp}a_1^{\sharp}v_1^{\sharp}u_2^{\sharp}}\dot{a}_2]
   \dashv ([v_2^{\sharp}]\dashv u_3\dot{a}_3v_3)
   =[\underline{u_1^{\sharp}a_1^{\sharp}v_1^{\sharp}u_2^{\sharp}}\dot{a}_2]
   \dashv [\underline{v_2^{\sharp}u_3^{\sharp}a_3^{\sharp}v_3^{\sharp}}] &\\
   =&[\underline{u_1^{\sharp}a_1^{\sharp}v_1^{\sharp}u_2^{\sharp}}]
   \vdash [\dot{a}_2\underline{v_2^{\sharp}u_3^{\sharp}a_3^{\sharp}v_3^{\sharp}}] =(u_1\dot{a}_1v_1\vdash [u_2^{\sharp}])\vdash [\dot{a}_2\underline{v_2^{\sharp}u_3^{\sharp}a_3^{\sharp}v_3^{\sharp}}] &\\
   = &u_1\dot{a}_1v_1\vdash ([u_2^{\sharp}]\vdash [\dot{a}_2\underline{v_2^{\sharp}u_3^{\sharp}a_3^{\sharp}v_3^{\sharp}}])
   =u_1\dot{a}_1v_1\vdash (u_2\dot{a}_2v_2\dashv u_3\dot{a}_3v_3).&
\end{align*}

The proof is completed.
\end{proof}

 Now we are ready to prove the main result of our article, in which we construct the free products of arbitrary digroups.

\begin{theorem}\label{nffreeproduct}
With Definition~\ref{defi-fp}, $(CW(D),\dot{e})$ is the free product of the family of indexed digroups~$\{(D_i',e)\mid i\in I\}$.  Moreover, the group part of~$(CW(D),\dot{e})$  is~$\{[\underline{u}]\mid u \mbox{ is a reduced word in } J^+  \}$, which is isomorphic to~$\PJ$, and the halo part of~$(CW(D),\dot{e})$  is~$\{u\dot{a}v\in CW(D)\mid \underline{u^{\sharp}a^{\sharp}v^{\sharp}}=e\}$.
\end{theorem}
\begin{proof}
We first show that~$(CW(D), \dot{e})$ is a digroup. By Lemma~\ref{pre-dise},  we know that~$CW(D)$ becomes a disemigroup. We proceed to show that~$(CW(D),\dot{e})$ is a digroup. Clearly~$\dot{e}$ is a bar-unit in the disemigroup~$CW(D)$.  Moreover, for every element~$u\dot{a}v\in CW(D)$, since $\Pi_{i\in I}^{\ast}J_i$ is a group,  by the construction of the free products of groups,  there exists a reduced word~$w\in J^+$  such that
$$
\underline{\underline{u^{\sharp}a^{\sharp}v^{\sharp}}w}=\underline{w\underline{u^{\sharp}a^{\sharp}v^{\sharp}}}=e.
$$
Thus we obtain
$$
u\dot{a}v\vdash [w]=[\underline{u^{\sharp}a^{\sharp}v^{\sharp}w^{\sharp}}]
=[\underline{\underline{u^{\sharp}a^{\sharp}v^{\sharp}}w}]=[e]=\dot{e}, \ \
[w]\dashv u\dot{a}v=[\underline{w^{\sharp}u^{\sharp}a^{\sharp}v^{\sharp}}]
=[\underline{w\underline{u^{\sharp}a^{\sharp}v^{\sharp}}}]=[e]=\dot{e}.
$$
Therefore, $(CW(D),\dot{e})$ is a digroup.

Now we show that~$(CW(D),\dot{e})$ is the free product of~$\{(D_i',e)\mid i\in I\}$.
For every $i\in I$, assume that $\varphi_i: (D_i',e)\rightarrow (CW(D),\dot{e})$ is the injective digroup homomorphism which is defined by $\varphi_i(x)=[x]=\dot{x}$ for every $x\in D_i'$.  Let~$(G,e')$ be an arbitrary digroup such that there exists a digroup homomorphism $\psi_i: (D_i',e)\rightarrow (G,e')$ for every~$i\in I$.  Recall that for every element~$d$ in~$D\setminus\{e\}$, the notation~$\sigma(d)$ means the index $i$ in~$I$ such that~$d$ lies in~$D_{i}'$.
Then we define a map $\eta$ from~$ (CW(D),\dot{e}) $ to~$(G,e')$ by the rule~$\eta(\dot{e})=e'$ and
$$
a_1...a_m\dot{a}b_1...b_n \mapsto  \psi_{\sigma(a_1)}(a_1)\vdash\cdots \vdash\psi_{\sigma(a_m)}(a_m)\vdash
\psi_{\sigma(a)}(a)\dashv\psi_{\sigma(b_1)}(b_1)\dashv\cdots\dashv \psi_{\sigma(b_n)}(b_n),$$
where~$\dot{e}\neq a_1...a_m\dot{a}b_1...b_n\in CW(D)$, $m\geq 0$ and~$n\geq 0$. If~$m=0$ (resp.~$n=0$), then~$\psi_{\sigma(a_1)}(a_1)\vdash\cdots \vdash\psi_{\sigma(a_m)}(a_m)\vdash$ (resp.~$\dashv\psi_{\sigma(b_1)}(b_1)\dashv\cdots\dashv \psi_{\sigma(b_n)}(b_n)$) does not appear. In particular, we have~$\eta(\dot{a})=\psi_{\sigma(a)}(a)$ for every~$a\in D\setminus\{e\}$.

Before showing that~$\eta$ is a digroup homomorphism, we first note that the map~$\eta$ is compatible with the extended reduction rule in the following sense: For all~$a_1,\dots, a_m, b_1,\dots, b_n$ in~$J$ and for every~$a\in D$,
assume that~$a_{t}\delta a_{t+1}=c$ in~$D_i'$ for some~$i\in I$,  where~$\delta\in\{\vdash, \dashv\}$. Then clearly we have
\begin{align*}
  &\eta(\dot{a}_1)\vdash \cdots \vdash \eta(\dot{a}_m)
\vdash\eta(\dot{a})\dashv\eta(\dot{b}_1)\dashv \cdots \dashv \eta(\dot{b}_{n})&\\
= &\eta(\dot{a}_1)\vdash \cdots \vdash \eta(\dot{a}_{t-1})
\vdash(\psi_i(\dot{a}_{t})\delta\psi_i(\dot{a}_{t+1})) \vdash  \cdots \vdash \eta(\dot{a}_m)
\vdash\eta(\dot{a})\dashv\eta(\dot{b}_1)\dashv \cdots \dashv \eta(\dot{b}_{n})&\\
= &\eta(\dot{a}_1)\vdash \cdots \vdash \eta(\dot{a}_{t-1})
\vdash(\psi_i(\dot{a}_{t}\delta \dot{a}_{t+1})) \vdash  \cdots \vdash \eta(\dot{a}_m)
\vdash\eta(\dot{a})\dashv\eta(\dot{b}_1)\dashv \cdots \dashv \eta(\dot{b}_{n})&\\
= &\eta(\dot{a}_1)\vdash \cdots \vdash \eta(\dot{a}_{t-1})
\vdash\eta(\dot{c}) \vdash  \cdots \vdash \eta(\dot{a}_m)
\vdash\eta(\dot{a})\dashv\eta(\dot{b}_1)\dashv \cdots \dashv \eta(\dot{b}_{n}).&
\end{align*}
Similarly, if~$a_m$ and $a $ lie in the same digroup, or  if~$b_1,a$ or  $b_t,b_{t+1}$ lie in the same digroup, then a similar result holds. Therefore, if~$a_1...a_m$ and~$b_1...b_n$ are words in~$J^{+}$ and if~$a\in D\setminus{J}$, then no matter whether~$\underline{a_1...a_m}\dot{a}\underline{b_1...b_n}$ is a good center-word or not,  we have
\begin{equation}\label{eta}
 \eta([\underline{a_1...a_m}\dot{a}\underline{b_1...b_n}])=\eta(\dot{a}_1)\vdash \cdots \vdash \eta(\dot{a}_m)
\vdash\eta(\dot{a})\dashv\eta(\dot{b}_1)\dashv \cdots \dashv \eta(\dot{b}_{n}).
\end{equation}
In particular, we deduce
\begin{equation}\label{change}
 \eta([\underline{a_1...a_m}])=\eta(\dot{a}_1)\dashv \cdots \dashv \eta(\dot{a}_m)
\end{equation}
for all~$a_1,...,a_m\in J$.

We also note that, for every~$c_1\in D_i'$, we have
 $$\eta([c_1^{\sharp}])=\psi_i(c_1^{\sharp}) =\psi_i(c_1\vdash e)=\psi_i(c_1)\vdash e'=\eta(\dot{c}_1)\vdash e'.$$
  By Lemma~\ref{Kinyon}, it follows that~$\psi_i(c_1\vdash e)$ lies in the group part of~$G'$, and in particular,  for every~$a_p\in J$,  the element~$\eta(\dot{a}_p)$ lies in
the group part of~$G'$. So by \eqref{change} and by Lemma~\ref{Kinyon}, we obtain
 $$\eta([\underline{a_1...a_m}])=\eta(\dot{a}_1)\dashv \cdots \dashv \eta(\dot{a}_m)
 =\eta(\dot{a}_1)\vdash \cdots \vdash \eta(\dot{a}_t)\dashv \eta(\dot{a}_{t+1}) \dashv \cdots \dashv \eta(\dot{a}_m)$$
 for every~$t$ such that~$1\leq t\leq m-1$.

Now we are ready to show that~$\eta$ is a digroup homomorphism.  For all~$u_1\dot{c}_1v_1$ and~$u_2\dot{c}_2v_2$ in~$CW(D)$, where~$u_1=a_1...a_m$, $v_1=b_1...b_n$, $u_2=d_1...d_p$ and~$v_2=f_1...f_q$,  if~$c_2\in  J$, then we have~$p=0$ and
 \begin{align*}
   \eta(u_1\dot{c}_1v_1\vdash u_2\dot{c}_2v_2)
   =&\eta([\underline{u^{\sharp}_{1}c_{1}^{\sharp}v_1^{\sharp}u_2^{\sharp}
   c_2^{\sharp}v_2^{\sharp}]})&\\
    =&\eta(\dot{a}_1)\dashv \cdots \dashv \eta(\dot{a}_m)
\dashv \eta([c_1^{\sharp}])\dashv \cdots \dashv \eta(\dot{b}_{n})
\dashv \eta(\dot{c}_2)\dashv \cdots \dashv \eta(\dot{f}_{q}) &\\
=&(\eta(\dot{a}_1)\vdash \cdots \vdash \eta(\dot{a}_m)
\vdash  \eta([c_1^{\sharp}])\dashv \cdots \dashv \eta(\dot{b}_{n})) \vdash (\eta(\dot{c}_2)\dashv \cdots \dashv \eta(\dot{f}_{q})) &\\
=&(\eta(\dot{a}_1)\vdash \cdots \vdash \eta(\dot{a}_m)
\vdash  (\eta(\dot{c}_1)\vdash e')\dashv \cdots \dashv \eta(\dot{b}_{n})) \vdash \eta(u_2\dot{c}_2v_2) &\\
=&(\eta(\dot{a}_1)\vdash \cdots \vdash \eta(\dot{a}_m)
\vdash \eta(\dot{c}_1)\dashv \cdots \dashv \eta(\dot{b}_{n}) )\vdash \eta(u_2\dot{c}_2v_2) &\\
 =& \eta(u_1\dot{c}_1v_1) \vdash \eta(u_2\dot{c}_2v_2);&
 \end{align*}
On the other hand, if~$c_2\in D\setminus J$, then we have
 \begin{align*}
   \eta(u_1\dot{c}_1v_1\vdash u_2\dot{c}_2v_2)
   =&\eta([\underline{u^{\sharp}_{1}c_{1}^{\sharp}v_1^{\sharp}u_2^{\sharp}}
   \dot{c}_2v_2^{\sharp}])&\\
   =&(\eta(\dot{a}_1)\vdash \cdots \vdash \eta(\dot{a}_m)
\vdash \eta([c_1^{\sharp}])\vdash\eta(\dot{b}_1)\vdash \cdots \vdash \eta(\dot{b}_{n}) ) \vdash \eta(u_2\dot{c}_2v_2)&\\
    =&(\eta(\dot{a}_1)\vdash \cdots \vdash \eta(\dot{a}_m)
\vdash(\eta(\dot{c}_1)\vdash e')\vdash\eta(\dot{b}_1)\vdash \cdots \vdash \eta(\dot{b}_{n}) ) \vdash \eta(u_2\dot{c}_2v_2)&\\
   =&(\eta(\dot{a}_1)\vdash \cdots \vdash \eta(\dot{a}_m)
\vdash\eta(\dot{c}_1)\dashv\eta(\dot{b}_1)\dashv \cdots \dashv \eta(\dot{b}_{n})) \vdash \eta(u_2\dot{c}_2v_2)&\\
 =& \eta(u_1\dot{c}_1v_1) \vdash \eta(u_2\dot{c}_2v_2).&
 \end{align*}

Similarly, we have~$\eta(u_1\dot{c}_1v_1\dashv u_2\dot{c}_2v_2) =\eta(u_1\dot{c}_1v_1)\dashv \eta(u_2\dot{c}_2v_2)$. It follows that~$\eta$ is a digroup homomorphism from~$(CW(D), \dot{e})$ to~$(G, e')$.

Since the set~$\{\dot{a}\mid a\in D\}$ generates~$(CW(D), \dot{e})$, it is clear that~$\eta$ is the unique digroup homomorphism  from~$(CW(D), \dot{e})$ to~$(G, e')$ such that~$\eta \varphi_i=\psi_i$ for every~$i\in I$.  Therefore, $(CW(D),\dot{e})$ is the free product of the family of  digroups~$\{(D_i',e)\mid i\in I\}$.

  Assume that~$u\dot{a}v$ lies in the group part of $(CW(D),\dot{e})$. Then by Lemma~\ref{Kinyon}, we have
$$u\dot{a}v=u\dot{a}v \vdash \dot{e}
=[\underline{u^{\sharp}a^{\sharp}v^{\sharp}e^{\sharp}}]
=[\underline{u^{\sharp}a^{\sharp}v^{\sharp}}],$$
 which holds if and only if~$u$ is the empty word and~$\dot{a}v$ lies in~$CW(D)$ satisfying~$av\in J^+$, namely, $u\dot{a}v=[av]$ and~$av$ is a reduced word in~$ J^+$. Now the isomorphism is clear by the constructions of~$\PJ$ and~$(CW(D), \dot{e})$.

Assume that~$u\dot{a}v$ lies in the halo part of $(CW(D),\dot{e})$. Then we have $$u\dot{a}v\vdash \dot{e}=[\underline{u^{\sharp}a^{\sharp}v^{\sharp}e^{\sharp}}]
=[\underline{u^{\sharp}a^{\sharp}v^{\sharp}}]=\dot{e}.$$ It follows immediately that we have~$\underline{u^{\sharp}a^{\sharp}v^{\sharp}}=e$. On the other hand, by the construction of~$(CW(D), \dot{e})$,   it is clear that  every~$u\dot{a}v\in CW(D)$ satisfying~$\underline{u^{\sharp}a^{\sharp}v^{\sharp}}=e$ lies in the halo part of~$(CW(D),\dot{e})$.
\end{proof}

Theorem~\ref{nffreeproduct} shows that the free product of digroups extends the notion of the free product of groups in the sense that, if every digroup $(D_i',e)$ is a group, $i\in I$, then~$(CW(D), \dot{e})$ is the free product of the groups~$\{D_i'\mid i\in I\}$.

Now we offer two easy examples of the free product of a family of digroups. By Remark~\ref{rema-dj}, we may always assume that the considered digroups are disjoint.
\begin{example}
Let~$(D_1,e_1)$ and~$(D_2,e_2)$ be two disjoint digroups such that the group part of~$(D_2,e_2)$ contains exactly one element. Then with the conventions and notations of Remark~\ref{rema-dj}, we have~$J_2=\{e\}$, $J=J_1$ and~$D=D_1'\cup D_2'$.
Moreover, the set of all good center-words over $D$ is
$$
CW(D)=\{u\dot{a}v \mid u,v\in (J_1\cup\{\varepsilon\})\setminus\{e\}, a\in D, \mbox{ and if } a\in D_1'  \mbox{ then } u=v=\varepsilon\}.
$$
With Definition~\ref{defi-fp}, $(CW(D), \dot{e})$ is the free product of~$(D_1',e)$ and~$(D_2',e)$, which is isomorphic to the free product of~$(D_1,e_1)$ and~$(D_2,e_2)$.
\end{example}

\begin{example}
Let~$\{(D_i,e_i)\mid i\in I\}$ be a family of disjoint digroups such that the group part of~$(D_i,e_i)$ contains exactly one element for every~$i\in I$.
Then with the conventions and notations of Remark~\ref{rema-dj}, we have~$J=J_i=\{e\}$ for every~$i\in I$ and~$D=\cup_{i\in I} D_i'$.
Moreover,  the set of all good center-words over $D$ is
$$CW(D)=\{\dot{a}\mid a\in D\}.$$
With Definition~\ref{defi-fp}, $(CW(D),\dot{e})$ is the free product of~$\{(D_i',e)\mid i\in I\}$, which is isomorphic to the free product of~$\{(D_i,e_i)\mid i\in I\}$.
\end{example}

\section*{Acknowledgments}   We thank very much the anonymous referee for an extremely careful reading of the original version of this article. The referee's insightful comments and suggestions help us greatly improve the exposition of the article.

\newcommand{\noopsort}[1]{}

\end{document}